\documentclass[11pt,reqno]{amsart}
\usepackage{amssymb,amscd,amsbsy,mathrsfs}
\newcommand{\lb}{\linebreak}

\renewcommand{\a}{\alpha}
\renewcommand{\b}{\beta}
\newcommand{\g}{\gamma}

\newcommand{\e}{\varepsilon}
\newcommand{\vk}{\varkappa}
\newcommand{\z}{\zeta}

\newcommand{\s}{\sigma}
\renewcommand{\t}{\tau}

\newcommand{\f}{\varphi}

\newcommand{\D}{\Delta}
\renewcommand{\L}{\Lambda}

\newcommand{\E}{{\mathscr E}}

\newcommand{\h}{{\mathscr H}}
\newcommand{\I}{{\mathscr I}}

\newcommand{\K}{{\mathscr K}}

\newcommand{\1}{{\bf 1}}

\newcommand{\T}{{\Bbb T}}

\newcommand{\R}{{\Bbb R}}
\newcommand{\Z}{{\Bbb Z}}

\newcommand{\0}{{\boldsymbol{0}}}

\newcommand{\bs}{\boldsymbol}

\newcommand{\m}{{\boldsymbol m}}
\newcommand{\bS}{{\boldsymbol S}}

\newcommand{\rf}[1]{(\ref{#1})}

\newcommand{\df}{\stackrel{\mathrm{def}}{=}}

\newcommand{\supp}{\operatorname{supp}}

\newcommand{\trace}{\operatorname{trace}}
\newcommand{\rank}{\operatorname{rank}}
\newcommand{\const}{\operatorname{const}}

\newcommand{\eeq}{\end{equation}}
\newcommand{\beq}{\begin{equation}}
\newcommand{\bay}{\begin{eqnarray}}
\newcommand{\ba}{\begin{align*}}
\newcommand{\ea}{\end{align*}}
\newcommand{\ey}{\end{eqnarray}}
\newcommand{\bey}{\begin{eqnarray*}}
\newcommand{\eey}{\end{eqnarray*}}

\newcommand{\be}{\infty}

\newcommand{\bl}{\blacksquare}

\newcommand{\Pf}{{\bf Proof. }}

\newcommand{\ov}{\overline}

\newtheorem{thm}{\hspace{\parindent}Theorem}[section]

\newtheorem{cor}[thm]{\hspace{\parindent}Corollary}
\newtheorem{lem}[thm]{\hspace{\parindent}Lemma}

\usepackage{amsmath}

\makeatletter
\def\upintkern@{\mkern-7mu\mathchoice{\mkern-3.5mu}{}{}{}}
\def\upintdots@{\mathchoice{\mkern-4mu\@cdots\mkern-4mu}%
 {{\cdotp}\mkern1.5mu{\cdotp}\mkern1.5mu{\cdotp}}%
 {{\cdotp}\mkern1mu{\cdotp}\mkern1mu{\cdotp}}%
 {{\cdotp}\mkern1mu{\cdotp}\mkern1mu{\cdotp}}}

\newcommand{\UpMultiIntegral}[1]{%
  \edef\ints@c{\noexpand\upintop
    \ifnum#1=\z@\noexpand\upintdots@\else\noexpand\upintkern@\fi
    \ifnum#1>\tw@\noexpand\upintop\noexpand\upintkern@\fi
    \ifnum#1>\thr@@\noexpand\upintop\noexpand\upintkern@\fi
    \noexpand\upintop
    \noexpand\ilimits@
  }%
  \futurelet\@let@token\ints@a
}
\makeatother

\DeclareFontFamily{OMX}{mdbch}{}
\DeclareFontShape{OMX}{mdbch}{m}{n}{ <->s * [0.8]  mdbchr7v }{}
\DeclareFontShape{OMX}{mdbch}{b}{n}{ <->s * [0.8]  mdbchb7v }{}
\DeclareFontShape{OMX}{mdbch}{bx}{n}{<->ssub * mdbch/b/n}{}

\DeclareSymbolFont{uplargesymbols}{OMX}{mdbch}{m}{n}
\SetSymbolFont{uplargesymbols}{bold}{OMX}{mdbch}{b}{n}
\DeclareMathSymbol{\upintop}{\mathop}{uplargesymbols}{82}
\DeclareMathSymbol{\upointop}{\mathop}{uplargesymbols}{"48}

\DeclareFontEncoding{MDB}{}{}
\DeclareFontFamily{MDB}{mdbch}{}
\DeclareFontShape{MDB}{mdbch}{m}{n}{ <->s * [0.8]  mdbchrmb }{}
\DeclareFontShape{MDB}{mdbch}{b}{n}{ <->s * [0.8]  mdbchbmb }{}
\DeclareFontShape{MDB}{mdbch}{bx}{n}{<->ssub * mdbch/b/n}{}
\DeclareFontSubstitution{MDB}{cmr}{m}{n}
\DeclareSymbolFont{mathdesignB}{MDB}{mdbch}{m}{n}%
\SetSymbolFont{mathdesignB}{bold}{MDB}{mdbch}{b}{n}%
\DeclareMathSymbol{\upintclockwise}{\mathop}{mathdesignB}{128}
\DeclareMathSymbol{\upointclockwise}{\mathop}{mathdesignB}{130}
\DeclareMathSymbol{\upointctrclockwise}{\mathop}{mathdesignB}{132}
\DeclareMathSymbol{\upoiint}{\mathop}{mathdesignB}{134}
\DeclareMathSymbol{\upoiiint}{\mathop}{mathdesignB}{136}

\makeatletter
\newcommand{\upint}{\DOTSI\upintop\ilimits@}
\newcommand{\upoint}{\DOTSI\upointop\ilimits@}
\makeatother

\theoremstyle{remark}

\newtheorem*{rem*}{Remark}

\newcommand\dg{\frak D}

\newcommand\Wp{\stackrel{1}{W}}
\newcommand\Wv{\stackrel{2}{W}}

\newcommand\mB{\mathcal{B}}

\newcommand{\Ba}{\big(B_{\be,1}^1\big)_+(\T^2)}
\newcommand\CA{{\rm C}_{\rm A}}

\begin{document}

\newcommand{\vse}{\vspace{.2in}}
\numberwithin{equation}{section}

\title{Functions of perturbed pairs of noncommuting contractions}
\author{A.B. Aleksandrov and V.V. Peller}
\thanks{The research of the first author is supposed by RFBR grant 17-01-00607.
The publication was prepared with the support of the
RUDN University Program 5-100}
\thanks{Corresponding author: V.V. Peller; email: peller@math.msu.edu}

\begin{abstract}
We consider functions $f(T,R)$ of pairs of noncommuting contractions on Hilbert space and study the problem for which functions $f$ we have Lipschitz type estimates in Schatten--von Neumann norms. We prove that if $f$ belongs to the Besov class $\Ba$ of nalytic functions in the bidisk, then
we have a Lipschitz type estimate for functions $f(T,R)$ of pairs of not necessarily commuting contractions $(T,R)$ in the Schatten--von Neumann norms $\bS_p$ for $p\in[1,2]$. On the other hand, we show that for functions in $\Ba$ there are no Lipschitz such type estimates for $p>2$ as well as in the operator norm.
\end{abstract}

\maketitle

{\bf
\footnotesize
\tableofcontents
\normalsize
}

\setcounter{section}{0}
\section{\bf Introduction}
\setcounter{equation}{0}
\label{In}

\medskip

The purpose of this paper is to 
study the behavior of functions $f(T,R)$ of (not necessarily commuting) contractions  $T$ and $R$ under perturbation. 
We are going to obtain Lipschitz type estimates  in the Sachatten--von Neumann norms $\bS_p$, $1\le p\le2$, for functions $f$ in the Besov class $\Ba$ of analytic functions.
Note that functions $f(T,R)$ of noncommuting contractions can be defined in terms of double operator integrals with respect to semi-spectral measures, see \S\:\ref{dvoitro} below.

This paper can be considered as a continuation of the results of \cite{Pe1}--\cite{Pe7}, \cite{AP1}--\cite{AP4}, \cite{AP6}, \cite{APPS}, \cite{NP}, \cite{ANP}, \cite{PS} and \cite{KPSS} for functions of perturbed self-adjoint operators, contractions, normal operators, dissipative operators, functions of collections of commuting operators and functions of collections of noncommuting operators.

Recall that a Lipschitz function $f$ on $\R$ does not have to be {\it operator Lipschitz}, i.e., the condition $|f(x)-f(y)|\le\const|x-y|$, $x,\,y\in\R$, does not imply that
$$
\|f(A)-f(B)\|\le\const\|A-B\|
$$
for arbitrary self-adjoint operators (bounded or unbounded, does not matter)
$A$ and $B$. This was first established in \cite{F}. 

It turned out that functions in the (homogeneous) Besov space $B_{\be,1}^1(\R)$ are operator Lipschitz; this was established in \cite{Pe1} and \cite{Pe3}
(see \cite{Pee} for detailed information about Besov classes). We refer the reader to the recent survey \cite{AP4} for detailed information on operator Lipschitz functions. In particular,  \cite{AP4} presents various sufficient conditions and necessary conditions for a function on $\R$ to be operator Lipschitz. It is well known that if $f$ is an operator Lipschitz function on $\R$, and $A$ and $B$ are self-adjoint operators such that the difference $A-B$ belongs to the Schatten--von Neumann class $\bS_p$, $1\le p<\be$, then $f(A)-f(B)\in\bS_p$ and 
$\|f(A)-f(B)\|_{\bS_p}\le\const\|A-B\|_{\bS_p}$. Moreover, the constant on the right does not depend on $p$. In particular, this is true for functions
$f$ in the Besov class $B_{\be,1}^1(\R)$, i.e.,
\bay
\label{LiShNOL}
\|f(A)-f(B)\|_{\bS_p}\le\const\|f\|_{B_{\be,1}^1}\|A-B\|_{\bS_p},\quad 1\le p\le\be.
\ey

However, it was discovered in \cite{AP1} (see also \cite{FN}) that the situation becomes quite different if we 
replace the class of Lipschitz functions with the class $\L_\a(\R)$  of
H\"older functions of order $\a$, $0<\a<1$. Namely, the inequality 
$|f(x)-f(y)|\le\const|x-y|^\a$, $x,\,y\in\R$, implies that 
$$
\|f(A)-f(B)\|\le\const\|A-B\|^\a
$$
for arbitrary self-adjoint operators $A$ and $B$. Moreover, 
it was shown in \cite{AP2}
that if $A-B\in\bS_p$, $p>1$, and $f\in\L_a(\R)$, then $f(A)-f(B)\in\bS_{p/\a}$ and
$$
\|f(A)-f(B)\|_{\bS_{p/a}}\le\const\|A-B\|_{\bS_p}^\a
$$
for arbitrary self-adjoint operators $A$ and $B$.

Analogs of the above results for functions of normal operators, functions of contractions, functions of dissipative operators and functions of commuting collections of self-adjoint operators were obtained in \cite {Pe2}, \cite{AP3}, \cite{APPS}, \cite{NP}.

Note that it was shown in \cite{PS} that for $p\in(1,\be)$, inequality \rf{LiShNOL}
holds for arbitrary {\it Lipschitz} (not necessarily operator Lipschitz) functions $f$ with constant on the right that depends on $p$.  An analog of this result for functions of commuting self-adjoint operators was obtained in \cite{KPSS}.

In \cite{ANP} similar problems were considered for functions of two noncommuting self-adjoint operators (such functions can be defined in terms of double operator integrals, see \cite{ANP}). It was shown in \cite{ANP} that for functions $f$ on 
$\R^2$ in the (homogeneous) Besov class $B_{\be,1}^1(\R^2)$ and for $p\in[1,2]$, the following Lipschitz type estimate holds:
$$
\|f(A_1,B_1)-f(A_2,B_2)\|_{\bS_p}\le\const
\max\big\{\|A_1-A_2\|_{\bS_p},\|B_1-B_2\|_{\bS_p}\big\}
$$
for arbitrary pairs $(A_1,B_1)$ and $(A_2,B_2)$
of (not necessarily commuting) self-adjoint operators.

However, it was shown in \cite{ANP} that for $p>2$ there is no such Lipschitz type
estimate in the $\bS_p$ norm as well as in the operator norm. Moreover, it follows from the construction given in \cite{ANP} that for $p\in(2,\be]$ and for 
positive numbers $\e,\,\s,\,M$, there exists a function $f$ in $L^\be(\R^2)$ with Fourier transform supported in $[-\s,\s]\times[-\s,\s]$ such that
$$
\max\big\{\|A_1-A_2\|_{\bS_p},\|B_1-B_2\|_{\bS_p}\big\}<\e
$$
while
$$
\|f(A_1,B_1)-f(A_2,B_2)\|_{\bS_p}>M.
$$
Here we use the notation $\|\cdot\|_{\bS_\be}$ for operator norm.

This implies that unlike in the case of commuting operators, there cannot be any H\"older type estimates in the norm of $\bS_p$, $p>2$, for H\"older functions $f$ of order $\a$. Moreover, for $p>2$, there cannot be any estimate for $\|f(A_1,B_1)-f(A_2,B_2)\|_{\bS_p}$ for functions in the Besov class 
$B_{\be,q}^s(\R)$
for any $q>0$ and $s>0$.

On the other hand, it was observed by the anonymous referee of \cite{ANP}
that unlike in the case of commuting self-adjoint operators, there is no Lipschitz type estimates for $\|f(A_1,B_1)-f(A_2,B_2)\|_{\bS_2}$ for {\it Lipschitz functions} $f$ on $\R^2$, see \cite{ANP}.

Finally, let us mention that in the case of functions of triples of noncommuting operators there are no such Lipschitz type estimates for functions in the Besov class 
$B_{\be,1}^1(\R^3)$ in the norm of $\bS_p$ for any $p\in[1,\be]$. This was 
established in \cite{Pe7}.

In \S\:\ref{dvoitro} we give an introduction to double and triple operator integrals and we define functions $f(T,R)$ of noncommuting contractions. We define the Haagerup and Haagerup-like tensor products of three copies of the disk-algebra
$\CA$ and we define triple operator integrals whose integrands belong to such tensor products.

Lipschitz type estimates in Schatten--von Neumann norm will be obtained in \S\:\ref{LiptysvN}. We show that for $p\in[1,2]$ and for a function $f$ on $\T^2$ in the analytic Besov space $\Ba$, the following Lipschitz type inequality holds:
$$
\big\|f(T_1,R_1)-f(T_0,R_0)\big\|_{\bS_p}\le\const
\max\big\{\|T_1-T_0\|_{\bS_p},\|R_1-R_0\|_{\bS_p}\big\}
$$
for arbitrary pairs $(T_0,T_1)$ and $(R_0,R_1)$ of contractions.
Recall that similar inequality was established in \cite{ANP} for functions of self-adjoint operators. However, to obtain this inequality for functions of contractions, we need new algebraic formulae. Moreover, to obtain this inequality for functions of contractions, we offer an approach that does not use triple operator integrals. To be more precise, we reduce the inequality to the case of analytic polynomials $f$ and we integrate over finite sets, in which case triple operator integrals become finite sums. We establish explicit representations of the operator differences $f(T_1,R_1)-f(T_0,R_0)$ for analytic polynomials $f$ in terms of finite sums of elementary tensors
which allows us to estimate the $\bS_p$ norms.

However, we still use triple operator integrals to obtain in \S\:\ref{pretrioi} explicit formulae for the operator differences for arbitrary functions $f$ in $\Ba$.

In \S\:\ref{differ} we study differentiability properties in Schatten--von Neumann norms of the function
$$
t\mapsto f\big(T(t),R(t)\big)
$$
for $f$ in $\Ba$ and  contractive valued functions $t\mapsto T(t)$
and $t\mapsto R(t)$. We obtain explicit formulae for the derivative in terms of triple operator integrals. Again, to prove the existence of the derivative, we do not need triple operator integrals.

As in the case of functions of pairs self-adjoint operators (see \cite{ANP}), there are no Lipschitz type estimates in the norm of $\bS_p$, $p>2$,  for functions of pairs of not necessarily commuting contractions $f(T,R)$, $f\in\Ba$. This will be established in \S\:\ref{counter}. Note that the construction differs from the construction in the case of self-adjoint operators given in \cite{ANP}.

In \S\:\ref{zada} we state some open problems and in \S\:\ref{Bes} we give an introduction to Besov classes on polydisks.

We use the notation $\m$ for normalized Lebesgue measure on the unit circle $\T$ and the notation $\m_2$ for normalized Lebesgue measure on $\T^2$.

For simplicity we assume that we deal with separable Hilbert spaces.

\

\section{\bf Besov classes of periodic functions}
\label{Bes}

\

In this section we give a brief introduction to Besov spaces on the torus.

To define Besov spaces on the torus $\T^d$, we consider an infinitely differentiable function $w$ on $\R$ such
that
$$
w\ge0,\quad\supp w\subset\left[\frac12,2\right],\quad\mbox{and} \quad w(s)=1-w\left(\frac s2\right)\quad\mbox{for}\quad s\in[1,2].
$$
Let $W_n$, $n\ge0$, be the trigonometric polynomials defined by
$$
W_n(\z)\df\sum_{j\in\Z^d}w\left(\frac{|j|}{2^n}\right)\z^j,\quad n\ge1,
\quad W_0(\z)\df\sum_{\{j:|j|\le1\}}\z^j,
$$
where 
$$
\z=(\z_1,\cdots,\z_d)\in\T^d,\quad j=(j_1,\cdots,j_d),\quad\mbox{and}\quad
|j|=\big(|j_1|^2+\cdots+|j_d|^2\big)^{1/2}.
$$
For a distribution $f$ on $\T^d$ we put
\bay
\label{f=sumfn}
\label{fnWn}
f_n=f*W_n,\quad n\ge0.
\ey
It is easy to see that
\bay
\label{fSigmafn}
f=\sum_{n\ge0}f_n;
\ey
the series converges in the sense of distributions.
We say that $f$ belongs the {\it Besov class} $B_{p,q}^s(\T^d)$, $s>0$, 
$1\le p,\,q\le\be$, if
\bay
\label{Bperf}
\big\{2^{ns}\|f_n\|_{L^p}\big\}_{n\ge0}\in\ell^q.
\ey

The analytic subspace $\big(B_{p,q}^s\big)_+(\T^d)$ of $B_{p,q}^s(\T^d)$ consists of functions $f$ in $B_{p,q}^s(\T^d)$ for which the Fourier coefficients 
$\widehat f(j_1,\cdots,j_d)$ satisfy the equalities: 
\bay
\label{anafunapo}
\widehat f(j_1,\cdots,j_d)=0\quad\mbox{whenever}\quad
\min_{1\le k\le d}j_k<0.
\ey

We refer the reader to \cite{Pee} for more detailed information about Besov spaces.

\

\section{\bf Double and triple operator integrals \\ with respect to semi-spectral measures}
\setcounter{equation}{0}
\label{dvoitro}

\

{\bf 3.1.~Double operator integrals.}
In this section we give a brief introduction to double and triple operator integrals with respect to semi-spectral measures. Double operator integrals with respect to 
spectral measures are expressions of the form
\bay
\label{dvooipolu}
\iint\Phi(x,y)\,dE_1(x)Q\,dE_2(y),
\ey
where $E_1$ and $E_2$ are spectral measures, $Q$ is a linear operator and
$\Phi$ is a bounded measurable function.
They appeared first in \cite{DK}. Later Birman and Solomyak developed in
\cite{BS1}--\cite{BS3} a beautiful theory of double operator integrals.

Double operator integrals with respect to semi-spectral measures were defined in 
\cite{Pe2}, see also \cite{AP4} (recall that the definition of a {\it semi-spectral measure} differs from the definition of a spectral measure by replacing the condition that it takes values in the set of orthogonal projections with the condition that it takes values in the set of nonnegative contractions, see \cite{AP4} for more detail).

For the double operator integral to make sense for an arbitrary bounded linear operator $T$, we have to impose an additional assumption on $\Phi$. The natural class of such functions $\Phi$ is called the class of {\it Schur multipliers}, see \cite{Pe1}. There are various characterizations of the class of Schur multipliers. In particular, $\Phi$ is a Schur multiplier if and only if it belongs to the Haagerup tensor product $L^\be(E_1)\otimes_{\rm h}\!L^\be(E_2)$ of $L^\be(E_1)$ and $L^\be(E_2)$, i.e.,
it admits a representation of the form
\bay
\label{tenzpre}
\Phi(x,y)=\sum_j\f_j(x)\psi_j(y),
\ey
where the $\f_j$ and $\psi_j$ satisfy the condition
\bay
\label{Haatepro}
\sum_j|\f_j|^2\in L^\be(E_1)\quad\mbox{and}\quad\sum_j|\psi_j|^2\in L^\be(E_2).
\ey
In this case
\bay
\label{fladlya doi}
\iint\Phi(x,y)\,dE_1(x)Q\,dE_2(y)=
\sum_j\Big(\int\f_j\,dE_1\Big)Q\Big(\int\psi_j\,dE_2\Big);
\ey
the series converges in the weak operator topology.
The right-hand side of this equality does not depend on the choice of a representation
of $\Phi$ in \rf{tenzpre}. 

One can also {\it consider double operator integrals of the form {\em\rf{dvooipolu}} in the case when $E_1$ and $E_2$ are semi-spectral measures}. In this case, as in the case of spectral measures, formula \rf{fladlya doi} still holds under the same assumption \rf{Haatepro}.

It is easy to see that if $\Phi$ belongs to the projective tensor product
$L^\be(E_1)\widehat\otimes L^\be(E_2)$ of $L^\be(E_1)$ and $L^\be(E_2)$, i.e.,
$\Phi$ admits a representation of the form \rf{tenzpre} with $\f_j$ and $\psi_j$
satisfying
$$
\sum_j\|\f_j\|_{L^\be(E_1)}\|\psi_j\|_{L^\be(E_2)}<\be,
$$
then $\Phi$ is a Schur multiplier and \rf{fladlya doi} holds.

\medskip

{\bf 3.2.~The semi-spectral measures of contractions.}
Recall that if $T$ is a contraction (i.e., $\|T\|\le1$) on a Hilbert space $\h$, then by the Sz.-Nagy dilation theorem
(see \cite{SNF}),  $T$ has a unitary dilation, i.e., there exist a Hilbert space $\K$ that contains $\h$ and a unitary operator $U$ on $\K$ such that
$$
T^n=P_\h U^n\big|\h,\quad n\ge0,
$$
where $P_\h$ is the orthogonal projection onto $\h$. 

Among all unitary dilations of $T$ one can always select a {\it minimal} unitary dilation (in a natural sense) and all minimal unitary dilations are isomorphic, 
see \cite{SNF}.

The existence of a unitary dilation allows us to construct the natural functional calculus $f\mapsto f(T)$ for functions $f$ in the disk-algebra $\CA$ defined by
$$
f(T)=P_\h f(U)\big|\h=P_\h\left(\int_\T f(\z)\,dE_U(\z)\right)\Big|\h,\quad f\in\CA.
$$
where $E_U$ is the spectral measure of $U$.

Consider the operator set function $\E_T$ defined on the Borel subsets of the unit circle $\T$ by
$$
\E_T(\D)=P_\h E_U(\D)\big|\h,\quad\D\subset\T.
$$
Then $\E_T$ is a {\it semi-spectral measure}. It can be shown that it does not depend on the choice of a unitary dilation.
The semi-spectral measure $\E_T$ is called the {\it semi-spectral measure} of 
$T$.

\medskip

{\bf 3.3.~Functions of noncommuting contractions.} Let $f$ be a function on
the torus $\T^2$ that belongs to the Haagerup tensor product
$\CA\!\otimes_{\rm h}\!\CA$, i.e., $f$ admits a representation of the form
$$
f(\z,\t)=\sum_j\f_j(\z)\psi_j(\t),\quad\z,~\t\in\T,
$$
where $\f_j$, $\psi_j$ are functions in $\CA$ such that
$$
\sup_{\z\in\T}\sum_j|\f_j(\z)|^2<\be
\quad\mbox{and}\quad
\sup_{\t\in\T}\sum_j|\psi_j(\t)|^2<\be.                                                                                                                                  
$$
For a pair $(T,R)$ of (not necessarily commuting contractions),
the operator $f(T,R)$ is defined as the double operator integral
$$
\iint_{\T\times\T}f(\z,\t)\,d\E_T(\z)\,d\E_R(\t)=
\iint_{\T\times\T}f(\z,\t)\,d\E_T(\z)I\,d\E_R(\t).
$$

Note that if $f\in\Ba$, then $f\in\CA\!\otimes_{\rm h}\!\CA$, and so we can take functions $f(T,R)$ of contractions for an arbitrary function $f$ in $\Ba$.
Indeed, if  $f$ is an analytic polynomial in two variables of degree at most $N$ in each variable, then we can represent $f$ in the form
$$
f(\z,\t)=\sum_{j=0}^N\z^j\left(\sum_{k=0}^N\widehat f(j,k)\t^k\right).
$$
Thus $f$ belongs to the projective tensor product $\CA\hat\otimes\CA$ and
\bay
\label{tenzprepol}
\|f\|_{\CA\hat\otimes\CA}\le\sum_{j=0}^N\sup_\t
\left|\sum_{k=0}^N\widehat f(j,k)\t^k\right|
\le(1+N)\|f\|_{L^\be}
\ey
It follows easily from \rf{Bperf} that every function $f$ of Besov class $\Ba$ belongs to  
$\CA\hat\otimes\CA$, and so the operator $f(T,R)$ is well defined. Clearly,
\bay
\label{opfunnekom}
f(T,R)=\sum_{n\ge0}\sum_{j=0}^{2^{n+1}}T^j
\left(\sum_{k=0}^{2^{n+1}}\widehat f_n(j,k)R^k\right),
\ey
where $f_n$ is the polynomial defined by \rf{fnWn}. It follows immediately from 
\rf{tenzprepol} and \rf{Bperf} that the series converges absolutely in the operator norm. Note that 
formula \rf{opfunnekom} can be used as a definition of the functions $f(T,R)$ of noncommuting contractions in the case when $f\in\Ba$.

\medskip

{\bf 3.4.~Triple operator integrals.~Haagerup tensor products.} There are several approaches to multiple operator integrals. Triple operator integrals are expressions of the form 
$$
W_\Phi\df\iiint\Phi(x,y,z)\,dE_1(x)X\,dE_2(y)Y\,dE_3(z),
$$
where $\Phi$ is a bounded measurable function, $E_1$, $E_2$ and $E_3$ are spectral measures, and $X$ and $Y$ are bounded linear operators on Hilbert space.

In \cite{Pe4} triple (and more general, multiple) operator integrals were defined
for functions $\Phi$ in the integral projective product 
\mbox{$L^\be(E_1)\otimes_{\rm i}L^\be(E_2)\otimes_{\rm i}L^\be(E_3)$}. For such
functions $\Phi$, the following Schatten--von Neumann properties hold:
$$
\left\|\iiint\Phi\,dE_1X\,dE_2Y\,dE_3\right\|_{\bS_r}
\le\|\Phi\|_{L^\be\otimes_{\rm i}L^\be\otimes_{\rm i}L^\be}
\|X\|_{\bS_p}\|Y\|_{\bS_q},\quad \frac1r=\frac1p+\frac1q,
$$
whenever $1/p+1/q\le1$. Later in \cite{JTT} triple (and multiple) operator integrals
were defined for functions $\Phi$ in the Haagerup tensor product
\mbox{$L^\be(E_1)\otimes_{\rm h}L^\be(E_2)\otimes_{\rm h}L^\be(E_3)$}.
However, it turns out that under the assumption 
$\Phi\in L^\be\!\otimes_{\rm h}\!L^\be\!\otimes_{\rm h}\!L^\be$,
the conditions $X\in\bS_p$ and $Y\in\bS_q$ imply that 
$\iiint\Phi\,dE_1X\,dE_2Y\,dE_3\in\bS_r$, $1/r=1/p+1/q$,
only under the conditions that $p\ge2$ and $q\ge2$, see \cite{AP5} (see also
\cite{ANP}). Moreover, the following inequality holds:
$$
\left\|\iiint\Phi\,dE_1X\,dE_2Y\,dE_3\right\|_{\bS_r}
\le\|\Phi\|_{L^\be\otimes_{\rm h}L^\be\otimes_{\rm h}L^\be}
\|X\|_{\bS_p}\|Y\|_{\bS_q},\quad \frac1r=\frac1p+\frac1q,
$$
whenever $p\ge2$ and $q\ge2$, see \cite{AP5}.

Note also that to obtain Lipschitz type estimates for functions of noncommuting self-adjoint operators in \cite{ANP}, we had to use triple operator integrals with integrands $\Phi$ that do not belong to the Haagerup tensor product  $L^\be\!\otimes_{\rm h}\!L^\be\!\otimes_{\rm h}\!L^\be$. That is why we had to introduce in \cite{ANP} Haagerup-like tensor products of the first kind and of the second kind.

In this paper we are going to use triple operator integrals with integrands being continuous functions on $\T^3$ that belong to Haagerup and Haagerup-like tensor products of three copies of the disk-algebra $\CA$. We briefly define such tensor products and discuss inequalities we are going to use in the next section.

\medskip

{\bf Definition 1.} 
We say that a continuous function $\Phi$ on $\T^3$ belongs to the {\it Haagerup tensor product} $\CA\!\otimes_{\rm h}\!\CA\!\otimes_{\rm h}\!\CA$ if
$\Phi$ admits a representation
\bay
\label{Haatepre}
\Phi(\z,\t,\vk)=\sum_{j,k\ge0}\a_j(\z)\b_{jk}(\t)\g_k(\vk),
\qquad\z,\,\t,\,\vk\in\T,
\ey
where $\a_j$, $\b_{jk}$ and $\g_k$ are functions in $\CA$ such that 
\bay
\label{noHaatepr}
\sup_{\z\in\T}\left(\sum_{j\ge0}|\a_j(\z)|^2\right)^{1/2}
\sup_{\t\in\T}\big\|\{\b_{jk}(\t)\}_{j,k\ge0}\big\|_{\mB}\,
\sup_{\vk\in\T}\left(\sum_{k\ge0}|\g_k(\vk)|^2\right)^{1/2}<\be.
\ey
Here $\|\cdot\|_\mB$ stands for the operator norm of a matrix (finite or infinite) on the space $\ell^2$ or on a finite-dimensional Euclidean space. By definition, the
norm of $\Phi$ in $\CA\!\otimes_{\rm h}\!\CA\!\otimes_{\rm h}\!\CA$ is the infimum of the left-hand side of \rf{noHaatepr} over all representations of $\Phi$ in the form of \rf{Haatepre}.

\medskip

Suppose that $\Phi\in\CA\!\otimes_{\rm h}\!\CA\!\otimes_{\rm h}\!\CA$
and both \rf{Haatepre} and \rf{noHaatepr}  hold. Let $T_1$, $T_2$ and $T_3$
be contractions with semi-spectral measures $\E_{T_1}$, $\E_{T_2}$
and $\E_{T_3}$. Then for bounded linear operators $X$ and $Y$, we can 
define the triple operator integral
\bay
\label{tropin}
W_\Phi=\iiint\Phi\,d\E_{T_1}X\,d\E_{T_2}Y\,d\E_{T_3}
\ey
as
\begin{align*}
W_\Phi&\df\sum_{j,k}
\Big(\int\a_j(\z)\,d\E_{T_1}(\z)\Big)X
\Big(\int\b_{jk}(\t)\,d\E_{T_2}(\t)\Big)Y
\Big(\int\g_k(\vk)\,d\E_{T_3}(\vk)\Big)\\[.2cm]
&=\sum_{j,k}\a_j(T_1)X\b_{jk}(T_2)Y\g_k(T_3).
\end{align*}
It is easy to verify that the series converges in the weak operator topology
if we consider partial sums over rectangles. It can be shown in the same way as in the case of triple operator integrals with respect to spectral measures that the sum on the right does not depend on the choice of a representation of $\Phi$ in the form of \rf{Haatepre}, see Theorem 3.1 of \cite{ANP}.

We are going to use Lemma 3.2 of \cite{AP5}.
Suppose that $\{Z_j\}_{j\ge0}$ is a sequence of bounded linear operators on Hilbert space such that
\bay
\label{nrvadlyaAj}
\left\|\sum_{j\ge0}Z_j^*Z_j\right\|^{1/2}\le M\quad\mbox{and}\quad
\left\|\sum_{j\ge0}Z_jZ_j^*\right\|^{1/2}\le M.
\ey
Let $Q$ be a bounded linear operator. Consider the row ${\rm R}_{\{Z_j\}}(Q)$ and
the column ${\rm C}_{\{Z_j\}}(Q)$ defined by
$$
{\rm R}_{\{Z_j\}}(Q)\df\big(Z_0Q\:Z_1Q\:Z_2Q\:\cdots\big)
$$
and
$$
{\rm C}_{\{Z_j\}}(Q)\df\left(\begin{matrix}QZ_0\\QZ_1\\QZ_2\\\vdots\end{matrix}\right).
$$
Then by Lemma 3.2 of \cite{AP5}, for $p\in[2,\be]$, the following inequalities hold:
\bay
\label{RCAjSp}
\big\|{\rm R}_{\{Z_j\}}(Q)\big\|_{\bS_p}\le M\|Q\|_{\bS_p}\quad\mbox{and}\quad
\big\|{\rm C}_{\{Z_j\}}(Q)\big\|_{\bS_p}\le M\|Q\|_{\bS_p}
\ey
whenever $Q\in\bS_p$.

It is easy to verify that under the above assumptions
\bay
\label{WPhi}
W_\Phi={\rm R}_{\{\a_j(T_1)\}}(X)\,B\,{\rm C}_{\{\g_j(T_3)\}}(Y),
\ey
where $B$ is the operator matrix $\{\b_{jk}(T_2)\}_{j,k\ge0}$.

\begin{lem}
\label{norBjk}
Under the above hypotheses, 
$$
\|B\|\le\sup_{\t\in\T}\big\|\{\b_{jk}(\t)\}_{j,k\ge0}\big\|_{\mB}.
$$
\end{lem}

\Pf Let $U$ be a unitary dilation
of the contraction $T_2$ on a Hilbert space $\K$, $\K\supset\h$. Clearly, we can consider the space $\ell^2(\h)$ as a subspace of $\ell^2(\K)$. It is easy to see that 
$$
\{\b_{jk}(T_2)\}_{j,k\ge0}=P_{\ell^2(\h)}\{\b_{jk}(U)\}_{j,k\ge0}\big|\ell^2(\h),
$$
where $P_{\ell^2(\h)}$ is the orthogonal projection onto $\ell^2(\h)$. The result follows from the inequality 
$\|\{\b_{jk}(U)\}_{j,k\ge0}\|\le\sup_{\t\in\T}\big\|\{\b_{jk}(\t)\}_{j,k\ge0}\big\|_{\mB}$,
which is a consequence of the spectral theorem. $\bl$

\medskip

It follows from Lemma 3.2 of \cite{Pe8} that under the above assumptions,
inequalities \rf{nrvadlyaAj} hold for $Z_j=\a_j(T_1)$, $j\ge0$, with 
$M=\sup_{\z\in\T}\left(\sum_{j\ge0}|\a_j(\z)|^2\right)^{1/2}$ and for 
$Z_j=\g_j(T_3)$, $j\ge0$, with 
$M=\sup_{\z\in\T}\left(\sum_{j\ge0}|\g_j(\z)|^2\right)^{1/2}$.
This together with Lemma \rf{norBjk} and inequalities \rf{RCAjSp}
implies that under the above assumptions,
\begin{multline}
\label{nervodlyapr}
\big\|{\rm R}_{\{\a_j(T_1)\}}(X)\,B\,{\rm C}_{\{\g_j(T_3)\}}(Y)\big\|_{\bS_r}
\\[.2cm]
\le
\sup_{\z\in\T}\left(\sum_{j\ge0}|\a_j(\z)|^2\right)^{1/2}
\sup_{\t\in\T}\big\|\{\b_{jk}(\t)\}_{j,k\ge0}\big\|_{\mB}\,
\sup_{\vk\in\T}\left(\sum_{k\ge0}|\g_k(\vk)|^2\right)^{1/2}
\end{multline}
whenever $p\ge2$, $q\ge2$ and $1/r=1/p+1/q$.

The following theorem is an analog of the corresponding result for triple operator integrals with respect to spectral measures, see \cite{AP5}. It follows immediately from \rf{nervodlyapr}.

\begin{thm}
\label{troipolume}
Let $T_1$, $T_2$ and $T_3$ be contractions, and let
$X\in\bS_p$ and $Y\in\bS_q$, $2\le p\le\be$, $2\le q\le\be$.
Suppose that 
$\Phi\in\CA\!\otimes_{\rm h}\!\CA\!\otimes_{\rm h}\!\CA$. Then 
$W_\Phi\in\bS_r$, $1/r=1/p+1/q$, and
\end{thm}
$$
\left\|\iiint\Phi\,d\E_{T_1}X\,d\E_{T_2}Y\,d\E_{T_3}\right\|_{\bS_r}
\le\|\Phi\|_{\CA\!\otimes_{\rm h}\!\CA\!\otimes_{\rm h}\!\CA}
\|X\|_{\bS_p}\|Y\|_{\bS_q}.
$$
Recall that by $\bS_\be$ we mean the class of bounded linear operators.

\medskip

{\bf 3.5. Haagerup-like tensor products.}
We define here Haagerup-like tensor products of disk-algebras by analogy with
Haagerup-like tensor products of $L^\be$ spaces, see \cite{ANP}.

\medskip

{\bf Definition 2.} 
A continuous function $\Phi$ on $\T^3$
is said to belong to the {\it Haagerup-like tensor product 
$\CA\!\otimes_{\rm h}\!\CA\!\otimes^{\rm h}\!\CA$ of the first kind} if
it admits a representation
\bay
\label{tensprpe}
\Phi(\z,\t,\vk)=\sum_{j,k\ge0}\a_j(\z)\b_{k}(\t)\g_{jk}(\vk),
\qquad\z,\,\t,\,\vk\in\T,
\ey
where $\a_j$, $\b_k$ and $\g_{jk}$ are functions in $\CA$ such that 
$$
\sup_{\z\in\T}\left(\sum_{j\ge0}|\a_j(\z)|^2\right)^{1/2}
\sup_{\t\in\T}\left(\sum_{k\ge0}|\b_k(\t)|^2\right)^{1/2}\,
\sup_{\vk\in\T}\big\|\{\g_{jk}(\vk)\}_{j,k\ge0}\big\|_{\mB}<\be.
$$

Clearly, $\Phi\in\CA\!\otimes_{\rm h}\!\CA\!\otimes^{\rm h}\!\CA$
if and only if the function
$$
(z_1,z_2,z_3)\mapsto\Phi(z_3,z_1,z_2)
$$
belongs to the Haagerup tensor product
$\CA\!\otimes_{\rm h}\!\CA\!\otimes_{\rm h}\!\CA$.

\medskip

Similarly, we can define the Haagerup-like tensor product
$\CA\!\otimes^{\rm h}\!\CA\!\otimes_{\rm h}\!\CA$ of the second kind.

\medskip

{\bf Definition 3.} 
A continuous function $\Phi$ on $\T^3$
is said to belong to the {\it Haagerup-like tensor product 
$\CA\!\otimes^{\rm h}\!\CA\!\otimes_{\rm h}\!\CA$ of the second kind} if
it admits a representation
\bay
\label{tensprvt}
\Phi(\z,\t,\vk)=\sum_{j,k\ge0}\a_{jk}(\z)\b_{j}(\t)\g_k(\vk),
\qquad\z,\,\t,\,\vk\in\T,
\ey
where $\a_{jk}$, $\b_j$ and $\g_k$ are functions in $\CA$ such that 
$$
\sup_{\z\in\T}\big\|\{\a_{jk}(\z)\}_{j,k\ge0}\big\|_{\mB}\,
\sup_{\t\in\T}\left(\sum_{j\ge0}|\b_j(\t)|^2\right)^{1/2}\,
\sup_{\vk\in\T}\left(\sum_{k\ge0}|\g_k(\vk)|^2\right)^{1/2}<\be.
$$

Let us first consider the situation when $\Phi$ is defined by \rf{tensprpe}
or by \rf{tensprvt} with summation over a finite set.
In this case triple operator integrals of the form \rf{tropin}
can be defined for arbitrary bounded linear operators $X$ and $Y$ and for arbitrary contractions $T_1$, $T_2$ and $T_3$.

Suppose that 
\bay
\label{Fitipa1}
\Phi(\z,\t,\vk)=\sum_{j\in F_1}\sum_{k\in F_2}\a_j(\z)\b_{k}(\t)\g_{jk}(\vk),
\qquad\z,\,\t,\,\vk\in\T,
\qquad\a_j,\:\b_k,\:\g_{jk}\in\CA,
\ey
where  $F_1$ and $F_2$ are finite sets. We put
\bay
\label{konsumpeti}
\iiint\Phi\,d\E_{T_1}X\,d\E_{T_2}Y\,d\E_{T_3}\df
\sum_{j\in F_1}\sum_{k\in F_2}\a_j(T_1)X\b_{k}(T_2)Y\g_{jk}(T_3).
\ey
Suppose now that 
\bay
\label{Fitipa2}
\Phi(\z,\t,\vk)=\sum_{j\in F_1}\sum_{k\in F_2}\a_{jk}(\z)\b_{j}(\t)\g_k(\vk),
\qquad\z,\,\t,\,\vk\in\T,
\qquad\a_{jk},\:\b_j,\:\g_k\in\CA,
\ey
where  $F_1$ and $F_2$ are finite sets.
Then we put
\bay
\label{konsumvtti}
\iiint\Phi\,d\E_{T_1}X\,d\E_{T_2}Y\,d\E_{T_3}\df
\sum_{j\in F_1}\sum_{k\in F_2}\a_{jk}(T_1)X\b_{j}(T_2)Y\g_k(T_3).
\ey

The following estimate is a very special case of Theorem \ref{pervto} below. However, we have stated it separately because its proof is elementary and does not require the definition of triple operator integrals with integrands in Haagerup-like tensor products.

\begin{thm}
\label{otsHaapevto}
Let $X$ and $Y$ be bounded linear operators and let $T_1$, $T_2$ and $T_3$ are contractions. Suppose that $F_1$ and $F_2$ are finite sets. The following statements hold:

{\em(i)} Let $\Phi$ be given by {\em\rf{Fitipa1}}.
Suppose that  
$q\ge2$ and $1/r\df1/p+1/q\in[1/2,1]$. If $X\in\bS_p$ and $Y\in\bS_q$, then
the sum on the right of {\em\rf{konsumpeti}} belongs to 
$\bS_r$ and
\begin{multline*}
\left\|\sum_{j\in F_1}\sum_{k\in F_2}\a_j(T_1)X\b_{k}(T_2)Y\g_{jk}(T_3)\in\bS_r\right\|_{\bS_r}\le\\[.2cm]
\;\;\;\sup_{\z\in\T}\left(\sum_{j\in F_1}|\a_j(\z)|^2\right)^{1/2}
\!\!\sup_{\t\in\T}\left(\sum_{k\in F_2}|\b_k(\t)|^2\right)^{\1/2}
\!\!\sup_{\vk\in\T}\big\|\big\{\g_{jk}(\vk)\big\}_{j\in F_1,k\in F_2}\big\|_\mB
\|X\|_{\bS_p}\|Y\|_{\bS_q}.
\end{multline*}

{\em(ii)} Let $\Phi$ be given by {\em\rf{Fitipa2}}.
Suppose that  
$q\ge2$ and $1/r\df1/p+1/q\in[1/2,1]$. If $X\in\bS_p$ and $Y\in\bS_q$, then
the sum on the right of {\em\rf{konsumvtti}} belongs to 
$\bS_r$ and
\begin{multline*}
\left\|\sum_{j\in F_1}\sum_{k\in F_2}\a_{jk}(T_1)X\b_{j}(T_2)Y\g_k(T_3)\in\bS_r\right\|_{\bS_r}\le\\[.2cm]
\sup_{\z\in\T}\big\|\big\{\a_{jk}(\z)\big\}_{j\in F_1,k\in F_2}\big\|_\mB
\sup_{\t\in\T}\left(\sum_{j\in F_1}|\b_j(\t)|^2\right)^{1/2}
\!\!\sup_{\vk\in\T}\left(\sum_{k\in F_2}|\g_k(\vk)|^2\right)^{\1/2}
\!\!\|X\|_{\bS_p}\|Y\|_{\bS_q}.
\end{multline*}
\end{thm}

\Pf Let us prove (i). The proof of (ii) is the same. We are going to use a duality argument. Suppose that $Q\in\bS_{r'}$ and $\|Q\|_{\bS_{r'}}\le1$,
$1/r+1/r'=1$. We have
\begin{multline*}
\sup_Q\left|
\trace\left(Q\sum_{j\in F_1}\sum_{k\in F_2}\a_{jk}(T_1)X\b_{j}(T_2)Y\g_k(T_3)
\right)\right|\\[.2cm]
=
\sup_Q\left|
\trace\left(\sum_{j\in F_1}\sum_{k\in F_2}\g_k(T_3)Q\a_{jk}(T_1)X\b_{j}(T_2)
\right)Y\right|\\[.2cm]
\le\|Y\|_{\bS_q}\sup_Q
\left\|\sum_{j\in F_1}\sum_{k\in F_2}
\g_k(T_3)Q\a_{jk}(T_1)X\b_{j}(T_2)\right\|_{\bS_{q'}}.
\end{multline*}
Th result follows now from \rf{WPhi} and \rf{nervodlyapr}. $\bl$

\medskip

{\bf 3.6. Triple operator integrals with integrands in Haagerup-like
tensor products.}
We define triple operator integrals with integrands in 
$\CA\!\otimes_{\rm h}\!\CA\!\otimes^{\rm h}\!\CA$ by  analogy with triple operator integrals with respect to spectral measures, see \cite{ANP} and \cite{AP5}.
Let $\Phi\in\CA\!\otimes_{\rm h}\!\CA\!\otimes^{\rm h}\!\CA$ and let $p\in[1,2]$.
Suppose that $T_1$, $T_2$ and $T_3$ are contractions. For 
an operator $X$ of class $\bS_p$ and for
a bounded linear operator $Y$, we define the triple operator integral
\bay
\label{WPhi1}
\Wp_\Phi\df
\iint\!\!\upint\Phi(\z,\t,\vk)\,d\E_{T_1}(\z)X\,d\E_{T_2}(\t)Y\,d\E_{T_3}(\vk)
\ey
as the following continuous linear functional on $\bS_{p'}$,
$1/p+1/p'=1$ (on the class of compact operators in the case $p=1$):
$$
Q\mapsto
\trace\left(\left(
\iiint
\Phi(\z,\t,\vk)\,dE_{T_2}(\t)Y\,dE_{T_3}(\vk)Q\,dE_{T_1}(\z)
\right)X\right).
$$
Note that the triple operator integral 
$\iiint\Phi(\z,\t,\vk)\,dE_{T_2}(\t)Y\,dE_{T_3}(\vk)Q\,dE_{T_1}(\z)$
is well defined as the integrand belongs to the Haagerup tensor product
$\CA\!\otimes_{\rm h}\!\CA\!\otimes_{\rm h}\!\CA$.

Again, we can define triple operator integrals with integrands in 
$\CA\!\otimes^{\rm h}\!\CA\!\otimes_{\rm h}\!\CA$ by analogy with the case of spectral measures, see \cite{ANP} and \cite{AP5}.
Let $\Phi\in\CA\!\otimes^{\rm h}\!\CA\!\otimes_{\rm h}\!\CA$ and let $T_1$, $T_2$ and $T_3$ be contractions. Suppose that $X$ is a bounded linear operator
and $Y\in\bS_p$, $1\le p\le2$.
The triple operator integral
\bay
\label{WPhi2}
\Wv_\Phi\df
\upint\!\!\!\iint\Phi(\z,\t,\vk)\,d\E_{T_1}(\z)X\,d\E_{T_2}(\t)Y\,d\E_{T_3}(\vk)
\ey
is defined as the continuous linear functional
$$
Q\mapsto
\trace\left(\left(
\iiint\Phi(\z,\t,\vk)\,dE_3(\vk)Q\,dE_1(\z)X\,dE_2(\t)
\right)Y\right)
$$
on $\bS_{p'}$ (on the class of compact operators if $p=1$).

As in the case of spectral measures (see \cite{AP5}), the following theorem can be proved:

\begin{thm}
\label{pervto}
Suppose that $T_1$, $T_2$ and $T_3$ are contractions,
and let $X\in\bS_p$ and $Y\in\bS_q$. The following statements hold:

{\em(1)} Let $\Phi\in\CA\!\otimes_{\rm h}\!\CA\!\otimes^{\rm h}\!\CA$.
Suppose 
$q\ge2$ and $1/r\df1/p+1/q\in[1/2,1]$. If $X\in\bS_p$ and $Y\in\bS_q$, then
the operator 
$\Wp_\Phi$ in {\em\rf{WPhi1}} belongs to $\bS_r$ and
\bay
\label{rpq}
\left\|\Wp_\Phi\right\|_{\bS_r}\le\|\Phi\|_{\CA\!\otimes_{\rm h}\,\!\CA\!\otimes^{\rm h}\CA}
\|X\|_{\bS_p}\|Y\|_{\bS_q};
\ey

{\em(2)} Let $\Phi\in\CA\!\otimes^{\rm h}\!\CA\!\otimes_{\rm h}\!\CA$.
Suppose that  $p\ge2$ and $1/r\df1/p+1/q\in[1/2,1]$. If $X\in\bS_p$ and $Y\in\bS_q$, then
the operator  
$\Wv_\Phi$ in {\em\rf{WPhi2}} belongs to $\bS_r$ and
$$
\Big\|\Wv_\Phi\Big\|_{\bS_r}\le\|\Phi\|_{\CA\!\otimes^{\rm h}\CA\!\otimes_{\rm h}\CA}
\|X\|_{\bS_p}\|Y\|_{\bS_q}.
$$
\end{thm}

\

\section{\bf Lipschitz type estimates in Schatten--von Neumann norms}
\setcounter{equation}{0}
\label{LiptysvN}

\

In this section we obtain Lipschitz type estimates in the Schatten--von Neumann classes $\bS_p$ for $p\in[1,2]$ for functions of contractions. To obtain such estimates, we are going to use an elementary approach and obtain elementary formulae that involve only finite sums. 

Later we will need explicit expressions for operator differences, which will be obtained in the next section  in terms of triple operator integrals. Such formulae will be used in \S\,\ref{differ} to obtain formulae for operator derivatives.

Suppose that $f$ is a function  that belongs to the Besov space 
$\big(B_{\be,1}^1\big)_+(\T^2)$ of analytic functions
(see \S\:\ref{Bes}). As we have observed in Subsection 3.3,
we can define functions $f(T,R)$ for (not necessarily commuting) contractions $T$ and $R$ on Hilbert space by formula \rf{opfunnekom}.

For a differentiable function $f$ on $\T$, we use the notation
$\dg f$ for the divided difference:
$$
(\dg f)(\z,\t)\df
\left\{\begin{array}{ll}
\displaystyle{\frac{f(\z)-f(\t)}{\z-\t}},&\z\ne\t\\[.4cm]f'(\z),&\z=\t,
\end{array}\right.
\qquad\z,\;\t\in\T.
$$
For a differentiable function $f$ on $\T^2$, we define the divided differences
$\dg^{[1]}f$ and $\dg^{[2]}f$ by
$$
\big(\dg^{[1]}f\big)(\z_1,\z_2,\t)\df
\left\{\begin{array}{ll}
\displaystyle{\frac{f(\z_1,\t)-f(\z_2,\t)}{\z_1-\z_2}},&\z_1\ne\z_2,\\[.4cm]
\displaystyle{\frac{\partial f}{\partial\z}\Big|_{\z=\z_1}},&\z_1=\z_2,
\end{array}\right.
\qquad\z_1,\;\z_2,\:\t\in\T,
$$
and
$$
\big(\dg^{[2]}f\big)(\z,\t_1,\t_2)\df
\left\{\begin{array}{ll}
\displaystyle{\frac{f(\z,\t_1)-f(\z,\t_2)}{\t_1-\t_2}},&\t_1\ne\t_2,\\[.4cm]
\displaystyle{\frac{\partial f}{\partial\t}\Big|_{\t=\t_1}},&\t_1=\t_2,
\end{array}\right.
\qquad\z,\;\t_1,\:\t_2\in\T.
$$

We need several elementary identities.

Let $\Pi_m$ be the set of $m$th roots of 1: 
$$
\Pi_m\df\{\xi\in\T:\xi^m=1\}
$$
and let
$$
\Upsilon_m(\z)\df\frac{\z^{m}-1}{m(\z-1)}=\frac1{m}\sum_{k=0}^{m-1}\z^k,\quad\z\in\T.
$$

The following elementary formulae are well known. We give proofs for completeness.

\begin{lem}
\label{136}
Let $f$ and $g$ be analytic polynomials in one variable of degree less than $m$. Then
$$
\int_\T f\ov g\,d\m=\frac1m\sum_{\xi\in\Pi_m}f(\xi)\ov{g(\xi)}.
$$
In particular,
$$
\int_\T |f|^2\,d\m=\frac1m\sum_{\xi\in\Pi_m}|f(\xi)|^2.
$$
\end{lem}

\Pf It suffices to consider the case where $f(z)=z^j$ and $g(z)=z^k$ with $0\le j, k<m$.
Then $-m<j-k<m$ and 
$$
\sum_{\xi\in\Pi_m}\xi^j\,\ov\xi^k=\left\{\begin{array}{ll}0,&j\ne k\\[.2cm]
m,&j=k.
\end{array}\right.\quad\bl
$$

\begin{cor}
\label{summodkv}
\bey
\sum_{\xi\in\Pi_m}|\Upsilon_m(\z\bar\xi)|^2=1,
\quad \z\in\T.
\eey
\end{cor}

In the same way we can obtain similar formulae for polynomials in several variables.
We need only the case of two variables.

\begin{lem}
\label{137}
Let $f$ and $g$ be polynomials in two variables of degree less than $m$
in each variable. Then
$$
\int_{\T^2} f\ov g\,d\m_2=\frac1{m^2}\sum_{\xi,\eta\in\Pi_m}f(\xi,\eta)\ov{g(\xi,\eta)}.
$$
In particular,
$$
\int_{\T^2} |f|^2\,d\m_2=\frac1{m^2}\sum_{\xi,\eta\in\Pi_m}|f(\xi,\eta)|^2.
$$
\end{lem}

\Pf It suffices to consider the case when $f(\z,\t)=\z^{j_1}\t^{j_2}$ and 
$g(\z,\t)=\z^{k_1}\t^{k_2}$ with 
$0\le j_1, j_2, k_1, k_2<m$.
Then $-m<j_1-k_1, j_2-k_2<m$ and 
\bay
\label{summirj1j2k1k2}
\sum_{\xi,\eta\in\Pi_m}\xi^{j_1}\eta^{j_2}\,\ov\xi^{k_1}\,\ov\eta^{k_2}=\left\{\begin{array}{ll}0,&(j_1,j_2)\ne(k_1,k_2)\\[.2cm]
m^2,&(j_1,j_2)=(k_1,k_2).
\end{array}\right.\quad\bl
\ey


Suppose now that $(T_0,R_0)$ and $(T_1,R_1)$ are pairs of not necessarily commuting contractions.

\begin{thm}
\label{flydlyaraz}
Let $f$ be an analytic polynomial in two variable of degree at most $m$
in each variable.
Then
\begin{align}
\label{perra}
f(T_1,R_1)-f(T_0,R_1)
=\sum_{\xi,\eta\in\Pi_m}\Upsilon_m(\ov\xi T_1)(T_1-T_0)\,\Upsilon_m(\ov\eta T_0)\,(\dg^{[1]}f)(\xi,\eta,R_1)
\end{align}
and
\begin{align}
\label{vtorra}
f(T_0,R_1)-f(T_0,R_0)
=\sum_{\xi,\eta\in\Pi_m}(\dg^{[2]}f)(T_0,\xi,\eta)\,\Upsilon_m(\ov\xi R_1)(R_1-R_0)\,\Upsilon_m(\ov\eta R_0).
\end{align}
\end{thm}

We are going to establish \rf{perra}.
The proof of \rf{vtorra} is similar.

We need the following lemma.

\begin{lem}
\label{odnape}
Let $\f$ be an analytic polynomial in one variable of degree at most $m$. Then
$$
\f(T_1)-\f(T_0)=
\sum_{\xi,\eta\in\Pi_m}\Upsilon_m(\ov\xi T_1)(T_1-T_0)\,\Upsilon_m(\ov\eta T_0)(\dg\f)(\xi,\eta).
$$
\end{lem}

{\bf Proof of the lemma.} 
Let $0\le j, j_0, k, k_0<m$. Then
$$
\sum_{\xi,\eta\in\Pi_m}(\ov\xi T_1)^{j_0}\,(\ov\eta T_0)^{k_0}\xi^j\eta^k=
\left\{\begin{array}{ll}\displaystyle{m^2T_1^jT_0^k},&(j_0,k_0)=(j,k),\\[.4cm]
0,&(j_0,k_0)\ne(j,k).
\end{array}\right.
$$
Thus, 
$$
\sum_{\xi,\eta\in\Pi_m}\Upsilon_m(\ov\xi T_1)\,\Upsilon_m(\ov\eta T_0)\xi^j\eta^k=T_1^jT_0^k
$$
if $0\le j, k<n$.
Hence,
$$
\sum_{\xi,\eta\in\Pi_m}\Upsilon_m(\ov\xi T_1)T_1\,\Upsilon_m(\ov\eta T_0)\xi^j\eta^k=T_1^{j+1}T_0^k
$$
and
$$
\sum_{\xi,\eta\in\Pi_m}\Upsilon_m(\ov\xi T_1)\,T_0\Upsilon_m(\ov\eta T_0)\xi^j\eta^k=T_1^{j}T_0^{k+1}.
$$
It follows that
$$
\sum_{\xi,\eta\in\Pi_m}\Upsilon_m(\ov\xi T_1)(T_1-T_0)\,\Upsilon_m(\ov\eta T_0)\xi^j\eta^k=T_1^j(T_1-T_0)T_0^k
$$
whenever $0\le j, k<m$.

Let $\f=\sum\limits_{s=0}^ma_sz^s$. It is easy to see that
$$
(\dg\f)(z,w)=\sum_{j,k\ge0, j+k<m}a_{j+k+1}z^jw^k.
$$
Hence,
\begin{multline*}
\sum_{\xi,\eta\in\Pi_m}\Upsilon_m(\ov\xi T_1)(T_1-T_0)\,\Upsilon_m(\ov\eta T_0)
(\dg\f)(\xi,\eta)\\
=\sum_{j,k\ge0, j+k<m}a_{j+k+1}\sum_{\xi,\eta\in\Pi_m}\Upsilon_m(\ov\xi T_1)(T_1-T_0)\,\Upsilon_m(\ov\eta T_0)\xi^j\eta^k\\
=\sum_{j,k\ge0, j+k<m}a_{j+k+1}T_1^j(T_1-T_0)T_0^k=\f(T_1)-\f(T_0).\quad\bl
\end{multline*}

\medskip 

{\bf Proof of Theorem \ref{flydlyaraz}.} Clearly, it suffices to prove \rf{perra}
in the case when $f(z_1,z_2)=\f(z_1)z_2^j$, where $\f$ is a polynomial of one
variable of degree at most $n$ and $0\le j\le m$.
Clearly, in this case
$$
f(T_1,R_1)-f(T_0,R_1)=\big(\f(T_1)-\f(T_0)\big)R_1^j.
$$
On the other hand,
$$
(\dg f^{[1]})(\xi,\eta,R_1)=(\dg\f)(\xi,\eta)R_1^j.
$$
Identity \rf{perra} follows now from Lemma \ref{odnape}. $\bl$

\medskip

For $K\in L^2(\T^2)$,  we denote by $\I_K$ the integral operator 
on $L^2(\T)$ with kernel function $K$, i.e.,
$$
(\I_K\f)(\z)=\int_\T K(\z,\t)\f(\t)\,d\m(\t),\quad \f\in L^2(\T).
$$
The following lemma allows us to evaluate the operator norm
$\|\I_K\|_{\mB(L^2)}$ of this operator for polynomials $K$ of degree less than $m$ in each variable in terms of the operator norms of the matrix
$\{K(\z,\eta)\}_{\z,\eta\in\Pi_m}$.

\begin{lem} 
\label{intopera}
Let $K$ be an analytic polynomial in two variables of degree less than $m$ in each variable.
Then 
$$
\|\{K(\xi,\eta)\}_{\xi,\eta\in\Pi_m}\|_{\mB}=m\|\I_K\|_{\mB(L^2)}.
$$
\end{lem}

\Pf It is easy to see that
\bey 
\|\I_K\|_{\mB(L^2)}=\sup_{\|\varphi\|_{L^2}\le1,\|\psi\|_{L^2}\le1}\left|\iint_{\T\times\T}K(\z,\t)\,\ov{\varphi(\z)\psi(\t)}\,d\m(\z)\,d\m(\t)\right|\\
=\sup_{\|\varphi\|_{L^2}\le1,\|\psi\|_{L^2}\le1}\left|\iint_{\T\times\T}K(\z,\t)\,\ov{\varphi_m(z)\psi_m(\t)}\,d\m(\z)\,d\m(\t)\right|,
\eey
where $\varphi_m(z)=\sum\limits_{k=0}^{m-1}\widehat\varphi(k)z^k$ and $\psi_m(z)=\sum\limits_{k=0}^{m-1}\widehat\psi(k)z^k$.
Hence,
$$
\|\I_K\|_{\mB(L^2)}=
\sup\left|\iint_{\T\times\T}K(\z,w)\,\ov{\varphi(z)\psi(w)}\,d\m(\z)\,d\m(\t)\right|,
$$
where the supremum is taken over all polynomials $\varphi$ and $\psi$ 
in one variable of degree less than $m$ and such that 
$\|\varphi\|_{L^2}\le1$, $\|\psi\|_{L^2}\le1$. Next, by Lemma \ref{137},
for arbitrary polynomials $\varphi$ and $\psi$ with $\deg\varphi<m$ and $\deg\psi<m$,
we have
$$
\iint_{\T\times\T}K(\z,\t)\,\ov{\varphi(z)\psi(w)}\,d\m(\z)\,d\m(\t)=\frac1{m^2}\sum_{\xi,\eta\in\Pi_m}K(\xi,\eta)\,\ov{\varphi(\xi)\psi(\eta)}.
$$
It remains to observe that by Lemma \ref{136},
$\|\varphi\|_{L^2}\le1$ if and only if $\sum\limits_{\xi\in\Pi_m}|\varphi(\xi)|^2\le m$
and the same is true for $\psi$.
 $\bl$

\begin{thm} 
\label{matrazra}
Let $g$ be a polynomial in one variable of degree at most $m$. Then
$$
\|\{(\dg g)(\xi,\eta)\}_{\xi,\eta\in\Pi_m}\|_{\mB}\le m\|g\|_{L^\be}.
$$
\end{thm}

\Pf The result follows from Lemma \ref{intopera} and the inequality
$$
\|\I_{\dg g}\|_{\mB(L^2)}\le\|g\|_{L^\be},
$$
which is a consequence of the fact that $\|\I_{\dg g}\|_{\mB(L^2)}$ is equal to the norm of the Hankel operator $H_{\bar g}$ on the Hardy class $H^2$, see
\cite{Pe5}, Ch. 1, Th. 1.10.
$\bl$

\begin{cor}
\label{sleu}
Let $f$ be a trigonometric polynomial of degree at most $m$ in each variable
and let $p\in[1,2]$.
Suppose that $T_1,\,R_1,T_0,\,R_0$ are contractions such that $T_1-T_0\in\bS_p$ and $R_1-R_0\in\bS_p$. Then 
$$
\|f(T_1,R_1)-f(T_0,R_0)\|_{\bS_p}\le2m\|f\|_{L^\be}
\max\big\{\|T_1-T_0\|_{\bS_p},\|R_1-R_0\|_{\bS_p}\big\}.
$$
\end{cor}

\Pf Let us estimate $\|f(T_1,R_1)-f(T_0,R_1)\|_{\bS_p}$. The norm
$\|f(T_0,R_1)-f(T_0,R_0)\|_{\bS_p}$
can be estimated in the same way. The result is a consequence of
formula \rf{perra}, Theorem \ref{otsHaapevto}, Theorem \ref{matrazra}
and Corollare \ref{summodkv}. $\bl$

\medskip

Corollary \ref{sleu} allows us to establish a Lipschitz type inequality for functions in
\lb$\big(B_{\be,1}^1\big)_+(\T^2)$.

\begin{thm}
\label{p2unit}
Let $1\le p\le2$ and let $f\in \big(B_{\be,1}^1\big)_+(\T^2)$. Suppose that $T_1,\,R_1,T_0,\,R_0$ are contractions such that $T_1-T_0\in\bS_p$ and $R_1-R_0\in\bS_p$. Then 
$$
\|f(T_1,R_1)-f(T_0,R_0)\|_{\bS_p}\le\const\|f\|_{B_{\be,1}^1}
\max\big\{\|T_1-T_0\|_{\bS_p},\|R_1-R_0\|_{\bS_p}\big\}.
$$
\end{thm}

\Pf Indeed, the result follows immediately from Corollary \ref{sleu} and inequality
\rf{Bperf}. $\bl$

\

\section{\bf A representation of operator differences\\ in terms
of triple operator integrals}
\setcounter{equation}{0}
\label{pretrioi}

\

In this section we obtain an explicit formula for the operator differences
$f(T_1,R_1)-f(T_0,R_0)$, $f\in\Ba$, in terms of triple operator integrals.

\begin{thm}
\label{razdraquasziHaa}
Let $f\in\Ba$. Then 
$$
\dg^{[1]}f\in\CA\!\otimes_{\rm h}\!\CA\!\otimes^{\rm h}\!\CA
\quad\mbox{and}\quad
\dg^{[2]}f\in\CA\!\otimes^{\rm h}\!\CA\!\otimes_{\rm h}\!\CA.
$$
\end{thm}

\begin{lem}
\label{predg1idg2}
Let $f$ be an analytic polynomial in two variables of degree at most $m$
in each variable.
Then
\bay
\label{dg1f}
\big(\dg^{[1]}f\big)(\z_1,\z_2,\t)=
\sum_{\xi,\eta\in\Pi_m}\Upsilon_m(\z_1\ov\xi)\,\Upsilon_m(\z_2\ov\eta)
\big(\dg^{[1]}f\big)(\xi,\eta,\t)
\ey
and
\bay
\label{dg2f}
\big(\dg^{[2]}f\big)(\z,\t_1,\t_2)=\sum_{\xi,\eta\in\Pi_m}
(\dg^{[2]}f)(\z,\xi,\eta)
\Upsilon_m(\t_1\ov\xi)\Upsilon_m(\t_2\ov\eta).
\ey
\end{lem}

\Pf Both formulae \rf{dg1f} and \rf{dg2f} can be verified straightforwardly. However, we deduce them from Theorem \ref{flydlyaraz}.

Formula \rf{dg1f} follows immediately from formula \rf{perra}
if we consider the special case when $T_0$, $T_1$ and $R_1$ 
are the operators on the one-dimensional space of multiplication 
by $\z_2$, $\z_1$ and $\t$. Similarly, formula \rf{dg2f}
follows immediately from formula \rf{vtorra}. $\bl$

\begin{cor}
\label{otsquasiHaanorm}
Under the hypotheses of Lemma {\em\ref{predg1idg2}},
$$
\big\|\dg^{[1]}f\big\|_{\CA\!\otimes_{\rm h}\!\CA\!\otimes^{\rm h}\CA}
\le m\|f\|_{L^\be}\quad\mbox{and}\quad
\big\|\dg^{[2]}f\big\|_{\CA\!\otimes^{\rm h}\CA\!\otimes_{\rm h}\!\CA}
\le m\|f\|_{L^\be}.
$$
\end{cor}

\Pf The result is a consequence of Lemma \ref{predg1idg2},
Theorem \ref{matrazra}, Corollary \ref{summodkv} and Definitions 2 and 3 in 
\S\:\ref{dvoitro}. $\bl$

\medskip

{\bf Proof of Theorem \ref{razdraquasziHaa}.} The result follows immediately from Corollary \ref{otsquasiHaanorm} and inequality \rf{Bperf}. $\bl$

\begin{thm}
\label{predraquasiHaa}
Let $p\in[1,2]$.
Suppose that $T_0,\,R_0,\,T_1,\,R_1$ are contractions such that
$T_1-T_0\in\bS_p$ and $R_1-R_0\in\bS_p$. Then for $f\in\Ba$,
the following formula holds:
\begin{align}
\label{osnfuni}
f(T_1,R_1)&-f(T_0,R_0)\nonumber\\[.2cm]
&=
\iint\!\!\upint\big(\dg^{[1]}f\big)(\z_1,\z_2,\t)
\,dE_{T_1}(\z_1)(T_1-T_0)\,dE_{T_2}(\z_2)\,dE_{R_1}(\t),\nonumber\\[.2cm]
&+\upint\!\!\!\iint\big(\dg^{[2]}f\big)(\z,\t_1,\t_2)
\,dE_{T_2}(\z)\,dE_{R_1}(\t_1)(R_1-R_0)\,dE_{R_2}(\t_2).
\end{align}
\end{thm}

\Pf Suppose first that $f$ is an analytic polynomial in two variables of degree at most $m$ in each variable. In this case equality \rf{osnfuni} is a consequence of
Theorem \ref{flydlyaraz}, Lemma \ref{predg1idg2}
and the definition of triple operator integrals given in Subsection 3.5.

In the general case we represent $f$ by the series \rf{f=sumfn} and apply 
\rf{osnfuni} to each $f_n$. The result follows from \rf{Bperf}. $\bl$

\

\section{\bf Differentiability properties}
\setcounter{equation}{0}
\label{differ}

\

In this section we study differentiability properties of the map
\bay
\label{fT1T0R1R0}
t\mapsto f\big(T(t),R(t)\big)
\ey
in the norm of $\bS_p$, $1\le p\le2$,  for functions
$t\mapsto T(t)$ and $t\mapsto R(t)$ that take contractive values and are differentiable in $\bS_p$.

We say that an operator-valued function $\Psi$ defined on an interval $J$ is {\it differentiable} in $\bS_p$ if $\Phi(s)-\Phi(t)\in\bS_p$ for any $s,\,t\in J$, and 
the limit
$$
\lim_{h\to\0}\frac1h\big(\Psi(t+h)-\Psi(t)\big)\df\Phi'(t)
$$
exists in the norm of $\bS_p$ for each $t$ in $J$.

\begin{thm}
Let $p\in[1,2]$ and let $f\in\Ba$. Suppose that
$t\mapsto T(t)$ and $t\mapsto R(t)$ are operator-valued functions on an interval
$J$ that take contractive values and are differentiable
in $\bS_p$. Then the function {\em\rf{fT1T0R1R0}} is differentiable on $J$ in 
$\bS_p$ and
\begin{align*}
\frac d{dt}&f\big(T(t),R(t)\big)\Big|_{t=s}
\\[.2cm]
&=
\iint\!\!\upint\big(\dg^{[1]}f\big)(\z_1,\z_2,\t)
\,dE_{T(s)}(\z_1)T'(s)\,dE_{T(s)}(\z_2)\,dE_{R(s)}(\t)\nonumber\\[.2cm]
&+\upint\!\!\!\iint\big(\dg^{[2]}f\big)(\z,\t_1,\t_2)
\,dE_{T(s)}(\z)\,dE_{R(s)}(\t_1)R'(s)\,dE_{R(s)}(\t_2),
\end{align*}
$s\in J$.
\end{thm}

\Pf As before, it suffices to prove the result in the case when 
$f$ is an analytic polynomial of degree at most $m$ in each variable.
Suppose that $f$ is such a polynomial. Put $F(t)\df f\big(T(t),R(t)\big)$.
We have
\begin{align*}
F(s+h)&-F(s)\\[.2cm]
&=\sum_{\xi,\eta\in\Pi_m}\Upsilon_m\big(\ov\xi T(s+h)\big)\big(T(s+h)-T(s)\big)\Upsilon_m\big(\ov\eta T(s)\big)\big(\dg^{[1]}f\big)\big(\xi,\eta,R(s+h)\big)\\[.2cm]
&+
\sum_{\xi,\eta\in\Pi_m}\big(\dg^{[2]}f\big)\big(T(s),\xi,\eta\big)
\Upsilon_m\big(\ov\xi R(s+h)\big)
\big(R(s+h)-R(s)\big)\Upsilon_m\big(\ov\eta R(s)\big).
\end{align*}
Clearly,
$$
\lim_{h\to0}\frac1h\big(T(s+h)-T(s)\big)=T'(s)\quad\mbox{and}\quad
\lim_{h\to0}\frac1h\big(R(s+h)-R(s)\big)=R'(s)
$$
in the norm of $\bS_p$. On the other hand, it is easy to see that
$$
\lim_{h\to0}\Upsilon_m\big(\ov\xi T(s+h)\big)=\Upsilon_m\big(\ov\xi T(s)\big),\quad
\lim_{h\to0}\big(\dg^{[1]}f\big)\big(\xi,\eta,R(s+h)\big)
=\big(\dg^{[1]}f\big)\big(\xi,\eta,R(s)\big)
$$
and
$$
\lim_{h\to0}\Upsilon_m\big(\ov\xi R(s+h)\big)=\Upsilon_m\big(\ov\xi R(s)\big)
$$
in the operator norm. Hence,
\begin{align*}
F'(s)
&=\sum_{\xi,\eta\in\Pi_m}\Upsilon_m\big(\ov\xi T(s)\big)T'(s)\Upsilon_m\big(\ov\eta T(s)\big)\big(\dg^{[1]}f\big)\big(\xi,\eta,R(s)\big)\\[.2cm]
&+
\sum_{\xi,\eta\in\Pi_m}\big(\dg^{[2]}f\big)\big(T(s),\xi,\eta\big)
\Upsilon_m\big(\ov\xi R(s)\big)
R'(s)\Upsilon_m\big(\ov\eta R(s)\big).
\end{align*}
It follows now from Lemma \ref{predg1idg2} 
and from the definition of triple operator integrals given in \S\:\ref{dvoitro}
that the right-hand side is equal to
\begin{multline*}
\iint\!\!\upint\big(\dg^{[1]}f\big)(\z_1,\z_2,\t)
\,dE_{T(s)}(\z_1)T'(s)\,dE_{T(s)}(\z_2)\,dE_{R(s)}(\t)\\[.2cm]
+\upint\!\!\!\iint\big(\dg^{[2]}f\big)(\z,\t_1,\t_2)
\,dE_{T(s)}(\z)\,dE_{R(s)}(\t_1)R'(s)\,dE_{R(s)}(\t_2)
\end{multline*}
which completes the proof. $\bl$

\

\

\section{\bf The case $\bs{p>2}$}
\setcounter{equation}{0}
\label{counter}

\

In this section we show that unlike in the case $p\in[1,2]$, there are no Lipschitz type estimates in the norm of $\bS_p$ in the case when $p>2$ for functions $f(T,R)$, $f\in\Ba$, of not noncommuting contractions. In particular, there are no such Lipschitz type estimates for functions $f\in\Ba$ in the operator norm. Moreover, we show that for $p>2$, such Lipschitz type estimates do not hold even for functions $f$ in $\Ba$ and for pairs of noncommuting {\it unitary operators}.

Recall that similar results were obtained in \cite{ANP} for functions of noncommuting self-adjoint operators. However, in this paper we use a different construction to obtain results for functions of unitary operators.

\begin{lem} 
\label{inter}
For each matrix $\{a_{\xi\, \eta}\}_{\xi, \eta\in\Pi_m}$,
there exists an analytic polynomial $f$ in two variables of degree at most $2m-2$ in each variable such that
$f(\xi,\eta)=a_{\xi\, \eta}$ for all $\xi, \eta\in\Pi_m$ and $\|f\|_{L^\be(\T^2)}\le \sup\limits_{\xi, \eta\in\Pi_m}|a_{\xi\, \eta}|$.
\end{lem}

\Pf Put 
$$
f(z,w)\df\sum_{\xi, \eta\in\Pi_m}a_{\xi\, \eta}\Upsilon_m^2(z\ov\xi)\Upsilon_m^2(w\ov\eta).
$$
Clearly, $f(\xi,\eta)=a_{\xi\, \eta}$ for all $\xi, \eta\in\Pi_m$ and
\begin{align*}
|f(z,w)|&\le\sup_{\xi, \eta\in\Pi_m}|a_{\xi\, \eta}|\sum_{\xi, \eta\in\Pi_m}|\Upsilon_m(z\ov\xi)|^2|\Upsilon_m(w\ov\eta)|^2\\[.2cm]
&=
\sup_{\xi, \eta\in\Pi_m}|a_{\xi\, \eta}|\sum_{\xi\in\Pi_m}|\Upsilon_m(z\ov\xi)|^2\sum_{\eta\in\Pi_m}|\Upsilon_m(w\ov\eta)|^2
=\sup_{\xi, \eta\in\Pi_m}|a_{\xi\, \eta}|
\end{align*}
by Corollary \ref{summodkv}. $\bl$

\begin{lem} 
\label{fU1U2V}
For each $m\in\Bbb N$, there exists an analytic polynomial $f$ in two variables of degree at most $4m-2$ in each variable, and unitary operators $U_1$, $U_2$
and $V$
such that 
$$
\|f(U_1,V)-f(U_2,V)\|_{\bS_p}>\pi^{-1}m^{\frac32-\frac1p} \|f\|_{L^\be(\T^2)}\|U_1-U_2\|_{\bS_p}
$$
for every $p>0$.
\end{lem}

\Pf  One can  select orthonormal bases 
$\{g_\xi\}_{\xi\in\Pi_m}$ and $\{h_\eta\}_{\eta\in\Pi_m}$  in an $m$-dimensional Hilbert space $\h$ such that
 $|(g_\xi,h_\eta)|=m^{-\frac12}$ for all $\xi,\eta\in\Pi_m$.
 Indeed, let $\h$ be the subspace of $L^2(\T)$ of analytic polynomials of degree
 less than $m$. 
 We can put $g_\xi\df\sqrt{m}\Upsilon_m(z\ov\xi)$ and
 $h_\eta=z^k$, where $\eta=e^{2\pi{\rm i}k/m}$, $0\le k\le m-1$.

Consider the rank one projections $\{P_\xi\}_{\xi\in\Pi_m}$ and $\{Q_\eta\}_{\eta\in\Pi_m}$ defined by
$P_\xi v = (v,g_\xi)g_\xi$, $\xi\in\Pi_m$, and $Q_\eta v = (v,h_\eta)h_\eta$, $\eta\in\Pi_m$. 
We define the unitary operators $U_1$, $U_2$, and $V$ by
$$
U_1=\sum_{\xi\in\Pi_m}\xi P_\xi,\quad U_2=e^{\frac{\pi{\rm i}}m}U_1 
\quad\mbox{and}\quad
V=\sum_{\eta\in\Pi_m}\eta Q_\eta.
$$

By Lemma \ref{inter}, there exists an analytic polynomial $f$ in two variables of degree at most $4m-2$ in each variable such that
$f(\xi,\eta)=\sqrt m(g_\xi,h_\eta)$ for all $\xi, \eta\in\Pi_m$,
$f(\xi,\eta)=0$ for all $\xi\in\Pi_{2m}\setminus \Pi_m$, $\eta\in\Pi_m$
and $\|f\|_{L^\be(\T^2)}=1$.
Clearly, $f(U_2,V)=\0$ and $f(U_1,V)=\sum\limits_{\xi,\eta\in\Pi_m}f(\xi,\eta)P_\xi Q_\eta$.
We have
$$
( f(U_1,V)h_\eta, g_\xi)=f(\xi,\eta)(h_\eta, g_\xi)=\frac1{\sqrt m}.
$$
Hence, $\rank f(U_1,V)=1$ and
$$
\|f(U_1,V)-f(U_2,V)\|_{\bS_p}=\|f(U_1,V)\|_{\bS_p}=\|f(U_1,V)\|_{\bS_2}=\sqrt m.
$$
It remains to observe that $\|U_1-U_2\|_{\bS_p}=
\big|1-e^{\frac{\pi{\rm i}}m}\big|m^{\frac1p}<\pi m^{\frac1p-1}$. $\bl$

\medskip

{\bf Remark.} If we replace the polynomial $f$ constructed in the proof of 
Lemma \ref{fU1U2V} with the polynomial $g$ defined by
$$
g(z_1,z_2)=z_1^{4m-2}z_2^{4m-2}f(z_1,z_2),
$$
it will obviously satisfy the same inequality:
\bay
\label{gU1U2V}
\|g(U_1,V)-g(U_2,V)\|_{\bS_p}>
\pi^{-1}m^{\frac32-\frac1p} \|g\|_{L^\be(\T^2)}\|U_1-U_2\|_{\bS_p}.
\ey

\medskip

It is easy to deduce from \rf{Bperf} that for such polynomials $g$
$$
c_1m\|g\|_{L^\be(\T^2)}\le\|g\|_{B^\be_{\be,1}}\le c_2m\|g\|_{L^\be(\T^2)}
$$
for some constants $c_1$ and $c_2$.

This together with \rf{gU1U2V} implies the following result:

\begin{thm} 
Let $M>0$ and $2<p\le\be$. Then there exist unitary operators $U_1$, $U_2$, $V$
and an analytic polynomial $f$ in two variables such that 
$$
\|f(U_1,V)-f(U_2,V)\|_{\bS_p}>M\|f\|_{B^1_{\infty,1}(\T^2)}\|U_1-U_2\|_{\bS_p}.
$$
\end{thm} 

\

\section{\bf Open problems}
\setcounter{equation}{0}
\label{zada}

\

In this section we state open problems for functions of noncommuting contractions.

\medskip

{\bf Functions of triples of contractions.} Recall that it was shown in \cite{Pe7} that for $f\in B_{\be,1}^1(\R)$, there are no Lipschitz type estimates in the norm of $\bS_p$ for any $p>0$ for functions $f(A,B,C)$ of triples of noncommuting self-adjoint operators. We conjecture that the same must be true in the case of functions of triples of not necessarily commuting contractions. Note that the construction given in \cite{Pe7} does not generalize to the case of functions of contractions.

\medskip

{\bf Lipschitz functions of noncommuting contractions.}  Recall that an unknown referee of \cite{ANP} observed that for {\it Lipschitz} functions $f$ on the real line
there are no Lipschitz type estimates for functions $f(A,B)$ of noncommuting self-adjoint operators in the Hilbert--Schmidt norm. The construction is given in \cite{ANP}. We conjecture that the same result must hold in the case of functions of noncommuting contractions.

\medskip

{\bf Lipschitz type estimates for $\bs{p>2}$ and H\"older type estimates.}
It follows from results of \cite{ANP} that in the case of functions of noncommuting self-adjoint operators for any $s>0$, $q>0$ and $p>2$, there exist pairs of self-adjoint operators $(A_0,A_1)$ and $(B_0,B_1)$ and a function $f$ in the homogeneous Besov space $B_{\be,q}^s(\R)$  such that
$\|f(A_1,B_1)-f(A_0,B_0)\|_{\bS_p}$ can be arbitrarily large while 
$\max\{\|A_1-A_0\|_{\bS_p},\|B_1-B_0\|_{\bS_p}\}$ can be arbitrarily small.
In particular, the condition $f\in B_{\be,q}^s(\R)$ does not imply any Lipschitz or H\"older type estimates in the norm of $\bS_p$, $p>2$, for any positive $s$ and $q$.

It is easy to see that in the case of contractions the situation is different: for any $q>0$ and $p\ge1$, there exists $s>0$ such that the condition $f\in B_{\be,q}^s$
guarantees a Lipschitz type estimate for functions of not necessarily commuting contractions in $\bS_p$. 

It would be interesting to find optimal conditions on $f$ that would guarantee Lipschitz or H\"older type estimates in $\bS_p$ for a given $p$.

\

\

\noindent
\noindent
\begin{tabular}{p{7cm}p{15cm}}
A.B. Aleksandrov & V.V. Peller \\
St.Petersburg Branch & Department of Mathematics \\
Steklov Institute of Mathematics  & Michigan State University \\
Fontanka 27, 191023 St.Petersburg & East Lansing, Michigan 48824\\
Russia&USA\\
&and\\
&Peoples' Friendship University\\
& of Russia (RUDN University)\\
&6 Miklukho-Maklaya St., Moscow,\\
& 117198, Russian Federation
\end{tabular}

\end{document}